\documentclass[11pt,a4paper]{article}
\usepackage[utf8]{inputenc}
\usepackage[T1]{fontenc}
\usepackage[english]{babel}
\usepackage{lmodern}
\usepackage{amsmath, amssymb, amsthm}
\usepackage{booktabs}
\usepackage{cite}
\usepackage{hyperref}
\usepackage{microtype}

\theoremstyle{plain}
\newtheorem{theorem}{Theorem}[section]

\theoremstyle{definition}
\newtheorem{definition}{Definition}[section]
\newtheorem{remark}{Remark}[section]

\DeclareMathOperator{\Arg}{Arg}
\DeclareMathOperator{\Log}{Log}
\DeclareMathOperator{\Cov}{Cov}
\DeclareMathOperator{\Var}{Var}

\title{\textbf{Operative Representation of Quasi-Arithmetic Means, Complex Branch Structure, and Covariant Commutator Defects in Iterated Logarithmic Operation Hierarchies}}
\author{Uwe Heiber\thanks{University of Erfurt, Germany. E-mail: \texttt{uwe.heiber@uni-erfurt.de}}}
\date{November 19, 2026}

\begin{document}
	
	\emergencystretch=1.5em
	
	\maketitle
	
	\begin{abstract}
		We study the bidirectional family of commutative binary operations $\text{Op}_m^{(a)}[x,y] = f_a^{(-m)}\big(f_a^{(m)}(x) + f_a^{(m)}(y)\big)$ ($m \in \mathbb{Z}$) generated by iterating the scaled logarithmic Cauchy isomorphism $f_a(z) = \frac{1}{a}\ln(z)$ ($a > 0$). First, we show that for all $m \ge -1$ the associated quasi-arithmetic Kolmogorov--Nagumo means $S_{\text{Op}_m}(\mathbf{x})$ admit an internal operative two-stage representation $S_{\text{Op}_m}(\mathbf{x}) = \text{Op}_{m+1}^{(a)}\big[\text{Op}_m^{(a)}[\mathbf{x}], \, N_{\text{Op}_{m+1}}(n)\big]$ with normalization factor $N_{\text{Op}_{m+1}}(n) = f_a^{(-m)}(1/n)$, and prove that this internal operative normalizability via a distributive adjacent level characterizes logarithmic generators among continuously conjugated group operations. Second, since the normalization factor $N_{\text{Op}_0}(n) < 0$ ($n \ge 2$) leaves the positive real domain $\mathbb{R}_{>0}$ required for further iterated logarithms at $m \le -2$, we examine the analytic continuation onto the Riemann surface of the logarithm, contrasting the inner branch-index invariance of $\text{Op}_{-1}^{(a)}$ with the associativity defect $\omega(x,z)$ generated by $\text{Op}_2^{(a)}$. Third, for arbitrary levels $m \in \mathbb{Z}$ and sample sizes $n \in \mathbb{N}$, we establish the general covariance pushforward theorem $\text{Sp}_{\text{Op}_{m+1}}(\mathbf{x}, \mathbf{y}) = f_a^{(-m)}\!\big(\Cov(f_a^{(m)}(\mathbf{x}), f_a^{(m)}(\mathbf{y}))\big)$ governing the commutator defect between mean formation and the adjacent operation.
	\end{abstract}
	
	\section{Conjugated Operations and Double-Fixed-Point Calibration}
	
	In classical algebra, addition ($\text{Op}_0 := Ad$, $+$) on $\mathbb{R}$ and multiplication ($\text{Op}_1 := Mu$, $\cdot$) on $\mathbb{R}_{>0}$ form two commutative group structures coupled via the distributive law \cite{bosch2020}. Since classical exponentiation $x^y$ is neither commutative nor associative, commutative extensions of the hierarchy of operations via iterated logarithmic conjugations were investigated early on, notably by Bennett \cite{bennett1915} and, in the context of projective and non-Diophantine arithmetics, by Czachor \cite{czachor2016}. Closely related logarithmic-exponential couplings arise in the theory of functional equations \cite{aczel2006, kuczma2009}, quasi-arithmetic means \cite{kolmogorov1930, nagumo1930}, idempotent analysis \cite{maslov1987}, and psychophysics \cite{fechner1860}.
	
	In the following, we consider the bidirectional chain of commutative binary operations parametrized by a scaling parameter $a \in \mathbb{R}_{>0}$,
	\begin{equation}
		\label{eq:hierarchie_kette}
		\begin{aligned}
			\dots \quad &\longrightarrow \quad \text{Op}_{-1}^{(a)} \quad \longrightarrow \quad \text{Op}_{0}^{(a)}\,(+) \\
			&\longrightarrow \quad \text{Op}_{1}^{(a)}\,(\cdot) \quad \longrightarrow \quad \text{Op}_{2}^{(a)} \quad \longrightarrow \quad \dots
		\end{aligned}
	\end{equation}
	on the commutative scalar fields $\mathbb{R}$ and $\mathbb{C}$ (in higher-dimensional non-commutative division algebras such as the quaternions $\mathbb{H}$, both the commutativity of multiplication and the homomorphism property $e^{x+y} = e^x e^y$ fail). To establish a terminological link to the lecture \cite{heiber2026}, we note that the descending log-sum-exp operation $\text{Op}_{-1}^{(a)}$ is occasionally referred to as \textit{Nexumtion} ($Ne$) and the ascending logarithmic product operation $\text{Op}_{2}^{(a)}$ as \textit{Athiktion} ($At$); for a standardized mathematical exposition, however, we primarily use the systematic level notation $\text{Op}_m^{(a)}$ ($m \in \mathbb{Z}$). Binary operations are consistently written with square brackets $\text{Op}_m^{(a)}[x,y]$, whereas in the generator function $f_a^{(m)}(z)$ the subscript denotes the parameter $a$ and the superscript $(m)$ denotes $m$-fold functional iteration.
	
	\textit{Structure of the paper:} Following a summary of the conjugation foundations and the special double-fixed-point calibration $a = \frac{\ln(2)}{2}$ in Section~1, Section~2 investigates the operative representation of quasi-arithmetic means together with its characterization theorem. Since the iterated real normalization factors leave the positive real axis for $m \le -2$ and $n \ge 2$, Section~3 develops the analytic continuation onto the Riemann surface of the logarithm. Finally, Section~4 proves a general covariance pushforward theorem for the commutator defects between mean formation and the adjacent operation.
	
	\subsection{Derivation of Adjacent Operations from Distributivity and Cauchy's Equation}
	
	Requiring the distributive law
	\begin{equation}
		\label{eq:distributivgesetz_allg}
		U\big[ O[x,y], \, O[x,z] \big] = O\big[ x, \, U[y,z] \big]
	\end{equation}
	between a subordinate operation $U$ and its directly superordinate operation $O$, together with the existence of a continuous group isomorphism $f$ connecting both levels via the conjugation relation \cite{aczel2006, kuczma2009}
	\begin{equation}
		\label{eq:homomorphie_grund}
		f\big( O[x,y] \big) = U\big[ f(x), \, f(y) \big] \quad \iff \quad O[x,y] = f^{-1}\!\Big( U\big[ f(x), \, f(y) \big] \Big),
	\end{equation}
	yields immediately at the reference level $(U, O) = (\text{Op}_0, \text{Op}_1) = (Ad, Mu)$ on $\mathbb{R}_{>0}$ the logarithmic Cauchy functional equation $f(x \cdot y) = f(x) + f(y)$. (More general deformed base operations $Ad_p[x,y] = (x^p + y^p)^{1/p}$ and $Mu_q[x,y] = qxy$ reduce to the classical pair by preserving the neutral element $Mu_q[x,1]=x \implies q=1$ and the diagonal condition $Ad_p[2,2]=4 \implies p=1$.) Excluding discontinuous Hamel solutions \cite{hamel1905}, differentiation $y f'(xy) = f'(x)$, evaluation at $x=1$ with $f'(1) = \frac{1}{a}$ ($a \neq 0$), and integration under the normalization condition $f(1) = 0$ show that the unique continuous family of solutions on $\mathbb{R}_{>0}$ is the scaled logarithm
	\begin{equation}
		\label{eq:generator_fa}
		f_a(z) = \frac{1}{a}\ln(z), \qquad f_a^{-1}(u) = e^{a u} \qquad (a \in \mathbb{R}_{>0}).
	\end{equation}
	Because the exponential function $u \mapsto e^{au}$ is globally injective only on the reals, the basic algebraic construction is carried out first on $\mathbb{R}$ and $\mathbb{R}_{>0}$; its analytic continuation to the complex covering surface is uniquely determined by the identity theorem of complex analysis and is examined in Section~3.
	
	\begin{theorem}[Representation of the adjacent levels $\text{Op}_{-1}^{(a)}$ and $\text{Op}_{2}^{(a)}$ from distributivity]
		\label{thm:eindeutigkeit}
		Let $f_a(z) = \frac{1}{a}\ln(z)$ ($a > 0$) be the continuous Cauchy isomorphism between $(\mathbb{R}_{>0}, \cdot)$ and $(\mathbb{R}, +)$.
		\begin{enumerate}
			\item The subordinate operation $\text{Op}_{-1}^{(a)}$ obtained by conjugating addition $\text{Op}_0[u,v] = u + v$ via $f_a^{-1}$ admits (for $e^{ax} + e^{ay} \neq 0$) the representation
			\begin{equation}
				\label{eq:ne_roh}
				\text{Op}_{-1}^{(a)}[x,y] = f_a\!\Big( \text{Op}_0\big[ f_a^{-1}(x), \, f_a^{-1}(y) \big] \Big) = \frac{1}{a}\ln\!\left( e^{a x} + e^{a y} \right).
			\end{equation}
			\item Let $\text{Op}_2^{(a)}: \mathbb{R}_{>0} \times \mathbb{R}_{>0} \to \mathbb{R}_{>0}$ be a continuous commutative operation of the form $\text{Op}_2^{(a)}[x,y] = f_a^{-1}\!\big(g(f_a(x), f_a(y))\big)$ with continuous $g: \mathbb{R}^2 \to \mathbb{R}$. If $\text{Op}_2^{(a)}$ is distributive over multiplication ($\text{Op}_2^{(a)}[x, y \cdot z] = \text{Op}_2^{(a)}[x,y] \cdot \text{Op}_2^{(a)}[x,z]$), then $g(u,v) = c \cdot u \cdot v$ ($c \in \mathbb{R}$). Under the normalization convention $g(1,1) = 1$, it follows that $c = 1$ and hence
			\begin{equation}
				\label{eq:at_roh}
				\text{Op}_2^{(a)}[x,y] = f_a^{-1}\!\big( f_a(x) \cdot f_a(y) \big) = \exp\!\left( \frac{1}{a}\ln(x)\ln(y) \right).
			\end{equation}
		\end{enumerate}
	\end{theorem}
	
	\begin{proof}
		Part~1 follows directly by substituting $f_a^{-1}(x) = e^{ax}$ into \eqref{eq:homomorphie_grund}. For Part~2, applying $f_a$ to $\text{Op}_2^{(a)}[x, y \cdot z] = \text{Op}_2^{(a)}[x,y] \cdot \text{Op}_2^{(a)}[x,z]$ with $u = f_a(x)$, $v = f_a(y)$, and $w = f_a(z)$ transforms the distributive law into the additive Cauchy equation $g(u, v+w) = g(u,v) + g(u,w)$ on $\mathbb{R}$. Continuity and symmetry $g(u,v) = g(v,u)$ imply $g(u,v) = \psi(u)v = c u v$ with $c = g(1,1) = 1$ \cite{aczel2006, kuczma2009}.
	\end{proof}
	
	\subsection{General Level Structure, Boundary Asymptotics, and Double-Fixed-Point Calibration}
	
	\begin{definition}[Bidirectional logarithmic hierarchy of operations]
		\label{def:isomorphismus_allg}
		For $a > 0$ and $m \in \mathbb{Z}$, starting from $\text{Op}_0^{(a)}[x,y] = x+y$, the operation $\text{Op}_m^{(a)}$ is defined recursively on its maximal real domain $D_m$ by $f_a\big(\text{Op}_m^{(a)}[x,y]\big) = \text{Op}_{m-1}^{(a)}\big[f_a(x), f_a(y)\big]$, and thus explicitly by
		\begin{equation}
			\label{eq:op_m_rekursion}
			\text{Op}_m^{(a)}[x,y] = f_a^{(-m)}\!\Big( f_a^{(m)}(x) + f_a^{(m)}(y) \Big) = f_a^{-(m-1)}\!\Big( f_a^{(m-1)}(x) \cdot f_a^{(m-1)}(y) \Big).
		\end{equation}
	\end{definition}
	
	Since $f_a^{(m)}: (D_m, \text{Op}_m^{(a)}) \to (\mathbb{R}, +)$ is a group isomorphism, the group axioms carry over directly by transport of structure from $(\mathbb{R}, +)$ and $(\mathbb{R}_{>0}\setminus\{1\}, \text{Op}_2^{(a)})$: on $\mathbb{R} \cup \{-\infty\}$, $\text{Op}_{-1}^{(a)}$ forms a commutative monoid with neutral element $n_{\text{Op}_{-1}} = -\infty$ and (on $\overline{\mathbb{R}}$) absorbing element $d_{\text{Op}_{-1}} = +\infty$ \cite{maslov1987}. On $G = \mathbb{R}_{>0}\setminus\{1\}$, $(G, \text{Op}_2^{(a)})$ forms an abelian group for every $a > 0$ with neutral element $n_{\text{Op}_2}^{(a)} = f_a^{-1}(1) = e^a$, inverse element $x^{-1}_{\text{Op}_2} = f_a^{-1}(1/f_a(x)) = \exp(a^2/\ln x)$, and global absorbing element $d_{\text{Op}_2} = f_a^{-1}(0) = 1$ on $\mathbb{R}_{>0}$.
	
	\begin{remark}[Boundary asymptotics of $\text{Op}_2^{(a)}$ on either side of the threshold $x = 1$]
		\label{rem:rand_asymptotik}
		Whereas for $m \in \{-1, 0, 1\}$ the neutral element $n_{\text{Op}_m}$ becomes the global absorbing element $d_{\text{Op}_{m+1}} = n_{\text{Op}_m}$ upon ascending to the next level ($-\infty \to 0 \to 1$), the absorbing element $d_{\text{Op}_2} = 1$ splits the half-line $\mathbb{R}_{>0}$ into the subintervals $(0,1)$ with $\ln(x) < 0$ and $(1,\infty)$ with $\ln(x) > 0$. Consequently, the limiting behavior of $\text{Op}_2^{(a)}[x,y] = \exp(\frac{1}{a}\ln x \ln y)$ at the boundaries $y \to 0^+$ and $y \to +\infty$ reverses upon crossing $x = 1$:
		\begin{equation}
			\label{eq:at_rand_asymptotik}
			\begin{aligned}
				\lim_{y \to 0^+} \text{Op}_2^{(a)}[x,y] &= \begin{cases} 0 & (x > 1), \\ 1 & (x = 1), \\ +\infty & (0 < x < 1), \end{cases} \\[1ex]
				\lim_{y \to +\infty} \text{Op}_2^{(a)}[x,y] &= \begin{cases} +\infty & (x > 1), \\ 1 & (x = 1), \\ 0 & (0 < x < 1). \end{cases}
			\end{aligned}
		\end{equation}
	\end{remark}
	
	\begin{remark}[Canonical double-fixed-point calibration $a = \frac{\ln(2)}{2}$]
		\label{rem:doppelfixpunkt}
		Although all structural theorems in this paper hold for arbitrary scaling parameters $a > 0$, the specific parameter value $a = \frac{\ln(2)}{2} = \frac{\ln(4)}{4} \approx 0.34657$ (corresponding to the logarithm base $b = e^a = \sqrt{2}$) enjoys a remarkable number-theoretic fixed-point property: since $(z_1, z_2) = (2, 4)$ is the unique integer solution pair to $z_1^{z_2} = z_2^{z_1}$ ($2^4 = 4^2 = 16$, so $\frac{\ln 2}{2} = \frac{\ln 4}{4}$), $z_1 = 2$ and $z_2 = 4$ are simultaneous fixed points of the isomorphism $f_a$ and its inverse ($f_a(2) = 2$, $f_a(4) = 4$, and $e^{2a} = 2$, $e^{4a} = 4$) precisely when $a = \frac{\ln(2)}{2}$. Consequently, the unique positive diagonal intersection $2+2 = 2\cdot 2 = 4$ of addition and multiplication transfers via $\text{Op}_m^{(a)}[2,2] = f_a^{(-m)}\big(f_a^{(m)}(2) + f_a^{(m)}(2)\big) = f_a^{(-m)}(4) = 4$ to all levels $m \in \mathbb{Z}$, and along the slice $x=2$ the operation $\text{Op}_2^{(a)}[2,y] = y^{\ln(2)/a} = y^2$ reproduces the classical squaring function (with $n_{\text{Op}_2} = \sqrt{2}$). In the sequel, we write $\text{Op}_m$ without a superscript whenever this specific calibration $a = \frac{\ln(2)}{2}$ is assumed.
	\end{remark}
	
	\section{Operative Representation of Quasi-Arithmetic Means}
	
	The central structural result in the real domain concerns the interplay between the family of operations $\text{Op}_m^{(a)}$ and quasi-arithmetic means in the sense of Kolmogorov and Nagumo \cite{kolmogorov1930, nagumo1930}. Whereas in a general quasi-arithmetic mean $M_\varphi(x_1, \dots, x_n) = \varphi^{-1}\!\big(\frac{1}{n}\sum_{i=1}^n \varphi(x_i)\big)$ the normalization $\frac{1}{n}$ is performed externally as a division in the image space $\mathbb{R}$, within the logarithmic hierarchy for all $m \ge -1$ the directly superordinate operation $\text{Op}_{m+1}^{(a)}$ acts as an internal normalization operator for the $n$-fold application of the base operation $\text{Op}_m^{(a)}$.
	
	\begin{definition}[Structure-specific normalization factor and hierarchical mean]
		\label{def:mittelwert_hierarchie}
		Let $D_m \subseteq \mathbb{R}$ be the domain of the associative base operation $U = \text{Op}_m^{(a)}$ and $O = \text{Op}_{m+1}^{(a)}$ be the directly superordinate operation.
		\begin{enumerate}
			\item The \textbf{structure-specific normalization factor} $N_{\text{Op}_{m+1}}(n)$ is the element that maps the $n$-fold composition of identical elements $x \in D_m$ under $\text{Op}_{m+1}^{(a)}$ back to $x$ (idempotency condition):
			\begin{equation}
				\label{eq:normierungsfaktor_def}
				\text{Op}_{m+1}^{(a)}\!\Big[ \underbrace{\text{Op}_m^{(a)}[x, \dots, x]}_{n\text{ times}}, \; N_{\text{Op}_{m+1}}(n) \Big] = x.
			\end{equation}
			\item For a vector $\mathbf{x} = (x_1, \dots, x_n) \in D_m^n$, the \textbf{hierarchical mean} $S_{\text{Op}_m}$ of level $m$ is defined by
			\begin{equation}
				\label{eq:mittelwert_allgemein}
				S_{\text{Op}_m}(x_1, \dots, x_n) := \text{Op}_{m+1}^{(a)}\!\Big[ \text{Op}_m^{(a)}[x_1, \dots, x_n], \; N_{\text{Op}_{m+1}}(n) \Big].
			\end{equation}
		\end{enumerate}
	\end{definition}
	
	\begin{theorem}[Explicit representation of the normalization factor and real domain boundary]
		\label{thm:master_normierung}
		Let $a > 0$ and $f_a(z) = \frac{1}{a}\ln(z)$. Then the normalization factors $N_{\text{Op}_{m+1}}(n)$ ($n \in \mathbb{N}$) satisfy:
		\begin{enumerate}
			\item \textbf{Real existence and closed form for $m \ge -1$:} For every hierarchy level $m \ge -1$, there exists a real normalization factor $N_{\text{Op}_{m+1}}(n) \in D_{m+1}$ independent of $x$, uniquely given by
			\begin{equation}
				\label{eq:master_normierung}
				N_{\text{Op}_{m+1}}(n) = f_a^{(-m)}\!\left( \frac{1}{n} \right) = f_a^{-(m+1)}\!\left( -\frac{1}{a}\ln(n) \right).
			\end{equation}
			With this factor, the operative two-stage representation \eqref{eq:mittelwert_allgemein} coincides identically with the quasi-arithmetic Kolmogorov--Nagumo mean generated by $\varphi_m = f_a^{(m)}$:
			\begin{equation}
				\label{eq:kn_darstellung}
				S_{\text{Op}_m}(x_1, \dots, x_n) = f_a^{(-m)}\!\left( \frac{1}{n}\sum_{i=1}^{n} f_a^{(m)}(x_i) \right).
			\end{equation}
			\item \textbf{Departure from the real domain for $m \le -2$:} For all sample sizes $n \ge 2$, level $m = -1$ (with $N_{\text{Op}_0}(n) = -\frac{1}{a}\ln(n) < 0$) is the lowest level of the hierarchy whose normalization factor exists in the reals. At level $m = -2$ ($U = \text{Op}_{-2}^{(a)}, \, O = \text{Op}_{-1}^{(a)}$), the value $-\frac{1}{a}\ln(n) < 0$ for $n \ge 2$ lies outside the positive real domain $\mathbb{R}_{>0}$ of $f_a$, so the idempotency equation \eqref{eq:normierungsfaktor_def} has no real solution $N_{\text{Op}_{-1}}(n) \in \mathbb{R}$, but is solved in the complex plane on the principal branch by
			\begin{equation}
				\label{eq:komplexer_normierer_m_minus_2}
				\begin{aligned}
					N_{\text{Op}_{-1}}(n) &= \frac{1}{a}\Log_0\!\left(-\frac{1}{a}\ln(n)\right) \\
					&= \frac{1}{a}\ln\!\left(\frac{1}{a}\ln(n)\right) + \frac{i\pi}{a} \;\in\; \mathbb{C}\setminus\mathbb{R} \qquad (n \ge 2).
				\end{aligned}
			\end{equation}
		\end{enumerate}
	\end{theorem}
	
	\begin{proof}
		\textbf{Part 1:} By \eqref{eq:op_m_rekursion}, the $n$-fold composition of identical elements satisfies $f_a^{(m)}\big(\text{Op}_m^{(a)}[x, \dots, x]\big) = n \cdot f_a^{(m)}(x)$, while at the adjacent level $f_a^{(m)}\big(\text{Op}_{m+1}^{(a)}[u, v]\big) = f_a^{(m)}(u) \cdot f_a^{(m)}(v)$. Applying $f_a^{(m)}$ to the idempotency equation \eqref{eq:normierungsfaktor_def} yields in the image space
		\begin{equation*}
			n \cdot f_a^{(m)}(x) \cdot f_a^{(m)}\!\big(N_{\text{Op}_{m+1}}(n)\big) = f_a^{(m)}(x) \quad \iff \quad f_a^{(m)}\!\big(N_{\text{Op}_{m+1}}(n)\big) = \frac{1}{n}.
		\end{equation*}
		Since $\frac{1}{n} > 0$ lies in the domain $\mathbb{R}_{>0}$ of $f_a^{(1)}$ and $f_a^{(-k)}(\mathbb{R}) = \mathbb{R}_{>0}$ for all $k \ge 1$, $N_{\text{Op}_{m+1}}(n) = f_a^{(-m)}(1/n)$ is well-defined and real for all $m \ge -1$. Substituting into \eqref{eq:mittelwert_allgemein} gives $f_a^{(m)}\big(S_{\text{Op}_m}(\mathbf{x})\big) = \big(\sum_{i=1}^n f_a^{(m)}(x_i)\big) \cdot \frac{1}{n}$, and applying $f_a^{(-m)}$ yields \eqref{eq:kn_darstellung}.
		
		\textbf{Part 2:} For $m = -2$, the condition $f_a^{(-2)}\big(N_{\text{Op}_{-1}}(n)\big) = \frac{1}{n} \iff f_a^{(-1)}\big(N_{\text{Op}_{-1}}(n)\big) = f_a(1/n) = -\frac{1}{a}\ln(n)$ requires $\exp\!\big(a \cdot N_{\text{Op}_{-1}}(n)\big) = -\frac{1}{a}\ln(n)$. Since the right-hand side is strictly negative for $n \ge 2$ and $a > 0$, whereas $\exp(a c) > 0$ for all $c \in \mathbb{R}$, no real solution exists. Over $\mathbb{C}$, the equation $\exp(a \cdot N_{\text{Op}_{-1}}(n)) = e^{i\pi}\frac{1}{a}\ln(n)$ is solved on the principal sheet by \eqref{eq:komplexer_normierer_m_minus_2}.
	\end{proof}
	
	For the levels $m \in \{-1, 0, 1, 2\}$, equation \eqref{eq:master_normierung} produces the explicit expressions:
	\begin{itemize}
		\item \textbf{Level $m = -1$ (Exponential mean):} $N_{\text{Op}_0}(n) = -\frac{1}{a}\ln(n)$ and $S_{\text{Op}_{-1}}(\mathbf{x}) = \text{Op}_{-1}^{(a)}[\mathbf{x}] - \frac{1}{a}\ln(n) = \frac{1}{a}\ln\!\big(\frac{1}{n}\sum_{i=1}^n e^{a x_i}\big)$.
		\item \textbf{Level $m = 0$ (Arithmetic mean):} $N_{\text{Op}_1}(n) = \frac{1}{n}$ and $S_{\text{Op}_0}(\mathbf{x}) = \frac{1}{n}\sum_{i=1}^n x_i$.
		\item \textbf{Level $m = 1$ (Geometric mean):} $N_{\text{Op}_2}(n) = e^{a/n}$ (or $2^{\frac{1}{2n}}$ for $a = \frac{\ln 2}{2}$), allowing the geometric mean to be represented purely within the commutative operation family:
		\begin{equation}
			\label{eq:s_mu_explizit}
			S_{\text{Op}_1}(x_1, \dots, x_n) = \text{Op}_2^{(a)}\!\left[ \prod_{i=1}^n x_i, \; e^{a/n} \right] = \left( \prod_{i=1}^n x_i \right)^{\!\frac{1}{n}}.
		\end{equation}
		\item \textbf{Level $m = 2$ (Logarithmic-geometric mean):} $N_{\text{Op}_3}(n) = e^{a e^{a/n}}$ and $S_{\text{Op}_2}(\mathbf{x}) = \exp\!\big(a \prod_{i=1}^n (\frac{1}{a}\ln x_i)^{1/n}\big)$.
	\end{itemize}
	
	In the classical theory of quasi-arithmetic means due to Kolmogorov \cite{kolmogorov1930}, Nagumo \cite{nagumo1930}, and Aczél \cite{aczel2006}, mean formation is viewed as an analytic evaluation $\varphi^{-1}(\frac{1}{n}\sum_{i=1}^n \varphi(x_i))$ in the linear image space $\mathbb{R}$, with division by $n$ carried out externally. From an algebraic standpoint, one may ask for which continuously conjugated group operations this normalization can be realized \textit{internally} by a distributive adjacent operation $O$ of the hierarchy itself. The following theorem provides a complete answer:
	
	\begin{theorem}[Characterization of operative mean representation via a conjugated adjacent level]
		\label{thm:charakterisierung_mittel}
		Let $\varphi: D \to \mathbb{R}$ and $\psi: D_+ \to \mathbb{R}$ (where $D_+ := \{x \in D \mid \varphi(x) > 0\}$) be continuous, strictly monotone bijections generating the abelian group operations
		\begin{equation*}
			U[x,y] := \varphi^{-1}\!\big(\varphi(x) + \varphi(y)\big), \qquad O[x,y] := \psi^{-1}\!\big(\psi(x) + \psi(y)\big)
		\end{equation*}
		on $D$ and $D_+$, respectively, such that $f := \psi \circ \varphi^{-1}: \mathbb{R}_{>0} \to \mathbb{R}$ is a continuous bijection satisfying $f(1) = 0$. Then:
		\begin{enumerate}
			\item \textbf{Necessary and sufficient criterion:} A quasi-arithmetic mean $M_\varphi(x_1, \dots, x_n) = \varphi^{-1}\!\big(\frac{1}{n}\sum_{i=1}^n \varphi(x_i)\big)$ admits an operative representation on $D_+^n$ for all $n \in \mathbb{N}$ with a sequence $N_O(n) \in D_+$ of the form
			\begin{equation}
				\label{eq:op_normierbarkeit_satz}
				M_\varphi(x_1, \dots, x_n) = O\Big[ U[x_1, \dots, x_n], \, N_O(n) \Big] \qquad (\forall \mathbf{x} \in D_+^n, \; \forall n \in \mathbb{N})
			\end{equation}
			if and only if the transition isomorphism $f = \psi \circ \varphi^{-1}$ satisfies the logarithmic Cauchy equation $f(u \cdot v) = f(u) + f(v)$ on $\mathbb{R}_{>0}$, i.e., $f(u) = \frac{1}{a}\ln(u)$ ($a \neq 0$).
			\item \textbf{Simultaneous bi-invariance (distinction from Hölder means):} Whereas the power means $M_p(\mathbf{x}) = (\frac{1}{n}\sum x_i^p)^{1/p}$ \cite{aczel2006} for $p \neq 0, 1$ are merely singly homogeneous with respect to multiplication, every mean $S_{\text{Op}_m}$ of the logarithmic hierarchy simultaneously satisfies shift invariance with respect to $\text{Op}_m^{(a)}$ and scale invariance with respect to $\text{Op}_{m+1}^{(a)}$:
			\begin{align}
				S_{\text{Op}_m}\!\Big( \text{Op}_{m}^{(a)}[\lambda, x_1], \dots, \text{Op}_{m}^{(a)}[\lambda, x_n] \Big) &= \text{Op}_{m}^{(a)}\!\Big[ \lambda, \, S_{\text{Op}_m}(x_1, \dots, x_n) \Big], \label{eq:invarianz_u} \\
				S_{\text{Op}_m}\!\Big( \text{Op}_{m+1}^{(a)}[c, x_1], \dots, \text{Op}_{m+1}^{(a)}[c, x_n] \Big) &= \text{Op}_{m+1}^{(a)}\!\Big[ c, \, S_{\text{Op}_m}(x_1, \dots, x_n) \Big]. \label{eq:invarianz_o}
			\end{align}
		\end{enumerate}
	\end{theorem}
	
	\begin{proof}
		\textbf{Part 1:} Sufficiency for $f(u) = \frac{1}{a}\ln(u)$ follows with $N_O(n) = \varphi^{-1}(1/n)$ as in Theorem~\ref{thm:master_normierung}. Conversely, suppose \eqref{eq:op_normierbarkeit_satz} holds for all $\mathbf{x} \in D_+^n$ and $n \in \mathbb{N}$. Setting $u_i := \varphi(x_i) > 0$, $S := \sum_{i=1}^n u_i = \varphi(U[\mathbf{x}]) \in \mathbb{R}_{>0}$, and $c(n) := \psi(N_O(n)) \in \mathbb{R}$, applying $\psi = f \circ \varphi$ to \eqref{eq:op_normierbarkeit_satz} and using the surjectivity of $\varphi$ onto $\mathbb{R}_{>0}$ gives
		\begin{equation}
			\label{eq:beweis_schritt1}
			f\!\left( \frac{S}{n} \right) = f(S) + c(n) \qquad (\forall S \in \mathbb{R}_{>0}, \; \forall n \in \mathbb{N}).
		\end{equation}
		Setting $S = n$ and using $f(1) = 0$ immediately yields $c(n) = -f(n)$. With $S = m \cdot u$ ($m \in \mathbb{N}, \, u > 0$), \eqref{eq:beweis_schritt1} similarly gives $f(u) = f(m u) - f(m)$, i.e., $f(m u) = f(u) + f(m)$. For any positive rational number $q = \frac{m}{n} \in \mathbb{Q}_{>0}$ and any $u \in \mathbb{R}_{>0}$, we obtain $f(q u) = f(\frac{m u}{n}) = f(m u) - f(n) = f(u) + f(m) - f(n)$. Evaluating at $u = 1$ gives $f(q) = f(m) - f(n)$ and hence $f(q u) = f(q) + f(u)$ on $\mathbb{Q}_{>0} \times \mathbb{R}_{>0}$. Since $\mathbb{Q}_{>0}$ is dense in $\mathbb{R}_{>0}$ and $f = \psi \circ \varphi^{-1}$ is continuous, passing to the limit $q_k \to v \in \mathbb{R}_{>0}$ establishes Cauchy's equation $f(v \cdot u) = f(v) + f(u)$ for all $u,v \in \mathbb{R}_{>0}$, so $f(u) = \frac{1}{a}\ln(u)$ \cite{aczel2006, kuczma2009}.
		
		\textbf{Part 2:} Relations \eqref{eq:invarianz_u} and \eqref{eq:invarianz_o} follow immediately upon applying $f_a^{(m)}$ from the linearity of the arithmetic mean, $\frac{1}{n}\sum_{i=1}^n \big(f_a^{(m)}(\lambda) + f_a^{(m)}(x_i)\big) = f_a^{(m)}(\lambda) + \frac{1}{n}\sum_{i=1}^n f_a^{(m)}(x_i)$ and $\frac{1}{n}\sum_{i=1}^n \big(f_a^{(m)}(c)\cdot f_a^{(m)}(x_i)\big) = f_a^{(m)}(c)\cdot \frac{1}{n}\sum_{i=1}^n f_a^{(m)}(x_i)$.
	\end{proof}
	
	\section{Analytic Continuation and Branch Coupling on the Riemann Surface}
	
	As shown by Theorem~\ref{thm:master_normierung} (Part~2), mean normalization for $m \le -2$ and $n \ge 2$ already leads out of $\mathbb{R}_{>0}$ into the complex plane $\mathbb{C}$ \cite{freitag2006}. Due to the $2\pi i$-periodicity of the exponential function, the analytic continuation of the operations $\text{Op}_m^{(a)}$ onto the Riemann surface of the logarithm \cite{forster1981} is governed by the branch logarithms
	\begin{equation}
		\label{eq:def_zweiglogarithmus}
		\Log_k(z) := \ln|z| + i\big( \Arg(z) + 2\pi k \big) \qquad (\Arg(z) \in (-\pi, \pi], \;\; k \in \mathbb{Z}).
	\end{equation}
	A direct comparison between the descending level $\text{Op}_{-1}^{(a)}$ and the ascending level $\text{Op}_{2}^{(a)}$ reveals a fundamental dichotomy in how the branch indices $k \in \mathbb{Z}$ behave under nested evaluations.
	
	\subsection{Complex Continuation of \texorpdfstring{$\text{Op}_{-1}^{(a)}$}{Op\_-1}: Singularity Divisor and Branch Invariance}
	
	On the complement of the zero locus of $W(x,y) := e^{ax} + e^{ay}$, the $k$-th branch of $\text{Op}_{-1}^{(a)}$ for $a > 0$ is defined by
	\begin{equation}
		\label{eq:nexumtion_komplex}
		\text{Op}_{-1, k}^{(a)}[x,y] := \frac{1}{a}\Log_k\!\big( e^{ax} + e^{ay} \big) = \frac{1}{a}\Big( \ln\big|e^{ax} + e^{ay}\big| + i\,\Arg\big(e^{ax} + e^{ay}\big) + 2\pi i k \Big).
	\end{equation}
	
	\begin{theorem}[Analytic properties and inner branch invariance of $\text{Op}_{-1, k}^{(a)}$]
		\label{thm:komplex_ne_eigenschaften}
		For every parameter $a > 0$ and $k \in \mathbb{Z}$, the operation $\text{Op}_{-1, k}^{(a)}$ on $\mathbb{C}^2$ possesses the following properties:
		\begin{enumerate}
			\item \textbf{Singularity divisor and differential translation invariance:} The logarithmic branch singularities form the countable family of parallel affine hyperplanes
			\begin{equation}
				\label{eq:ne_singularitaet}
				\mathcal{S}_{-1} = \left\{ (x,y) \in \mathbb{C}^2 \;\middle|\; x - y = \frac{i\pi}{a}(2n + 1), \;\; n \in \mathbb{Z} \right\}.
			\end{equation}
			Outside $\mathcal{S}_{-1}$, the partial derivatives $\frac{\partial \text{Op}_{-1,k}^{(a)}}{\partial x} = \frac{e^{ax}}{e^{ax}+e^{ay}}$ and $\frac{\partial \text{Op}_{-1,k}^{(a)}}{\partial y} = \frac{e^{ay}}{e^{ax}+e^{ay}}$ are independent of the branch index $k$ and satisfy identically $\frac{\partial \text{Op}_{-1,k}^{(a)}}{\partial x} + \frac{\partial \text{Op}_{-1,k}^{(a)}}{\partial y} = 1$. A single positive loop around any hyperplane in $\mathcal{S}_{-1}$ induces the monodromy transition $\text{Op}_{-1, k}^{(a)} \mapsto \text{Op}_{-1, k+1}^{(a)} = \text{Op}_{-1, k}^{(a)} + \frac{2\pi i}{a}$, while the zero locus lies entirely on the principal sheet $k=0$ along $e^{ax} + e^{ay} = 1$.
			\item \textbf{Vanishing of the inner associativity defect:} Under nested evaluation for $(x,y,z) \in \mathbb{C}^3$ (away from the singular sets), for any fixed outer branch index $k_{\text{ext}} \in \mathbb{Z}$, the total value does \textit{not} depend on the inner branch indices $k_{xy}, k_{yz} \in \mathbb{Z}$:
			\begin{equation}
				\label{eq:ne_assoziativ_invarianz}
				\begin{aligned}
					\text{Op}_{-1, k_{\text{ext}}}^{(a)}\!\Big[ \text{Op}_{-1, k_{xy}}^{(a)}[x,y], \, z \Big] &= \text{Op}_{-1, k_{\text{ext}}}^{(a)}\!\Big[ x, \, \text{Op}_{-1, k_{yz}}^{(a)}[y,z] \Big] \\
					&= \frac{1}{a}\Log_{k_{\text{ext}}}\!\big( e^{ax} + e^{ay} + e^{az} \big).
				\end{aligned}
			\end{equation}
		\end{enumerate}
	\end{theorem}
	
	\begin{proof}
		Part~1 follows from $e^{ax}+e^{ay} = e^{ay}(e^{a(x-y)}+1) = 0 \iff a(x-y) = i\pi(2n+1)$ and elementary differentiation of \eqref{eq:nexumtion_komplex}. For Part~2, substituting the inner branch into the outer exponential yields
		\begin{align*}
			\exp\!\Big( a \cdot \text{Op}_{-1, k_{xy}}^{(a)}[x,y] \Big) &= \exp\!\Big( \Log_{k_{xy}}(e^{ax} + e^{ay}) \Big) \\
			&= (e^{ax} + e^{ay}) \cdot e^{2\pi i k_{xy}} = e^{ax} + e^{ay},
		\end{align*}
		so the inner branch indices $k_{xy}$ and $k_{yz}$ are annihilated by the $2\pi i$-periodicity of the exponential function.
	\end{proof}
	
	\subsection{Complex Continuation of \texorpdfstring{$\text{Op}_{2}^{(a)}$}{Op\_2}: Index Balance and Associator Structure}
	
	In sharp contrast to Theorem~\ref{thm:komplex_ne_eigenschaften} (Part~2), in the complex continuation of the ascending operation $\text{Op}_2^{(a)}$ on $S = \mathbb{C}\setminus\{0\}$,
	\begin{equation}
		\label{eq:athiktion_def_komplex}
		\text{Op}_{2; \, k_x, k_y}^{(a)}[x,y] := \exp\!\left( \frac{1}{a}\Log_{k_x}(x)\Log_{k_y}(y) \right) \qquad (x,y \in S, \;\; k_x, k_y \in \mathbb{Z}),
	\end{equation}
	the branch logarithms appear multiplicatively coupled in the exponent. While the complex distributive law $\text{Op}_{2; \, k_x, k_{yz}}^{(a)}[x, y \cdot z] = \text{Op}_{2; \, k_x, k_y}^{(a)}[x,y] \cdot \text{Op}_{2; \, k_x, k_z}^{(a)}[x,z]$ is preserved via $\Log_{k_y}(y) + \Log_{k_z}(z) = \Log_{k_y + k_z + m}(y \cdot z)$ through the linear index balance $k_{yz} = k_y + k_z + m$ with principal-value jump $m = \frac{1}{2\pi}\big(\Arg(y) + \Arg(z) - \Arg(y\cdot z)\big) \in \{-1, 0, 1\}$, nested evaluation gives rise to a non-trivial associativity defect:
	
	\begin{theorem}[Conditional complex associativity, factorization, and kernel of the associator $\omega(x,y,z)$]
		\label{thm:komplex_assoziativ}
		Let $x,y,z \in S = \mathbb{C}\setminus\{0\}$, $a > 0$, and $k_x, k_y, k_z, k_{xy}, k_{yz} \in \mathbb{Z}$ be fixed branch indices. Then the multiplicative \textbf{associativity defect} (associator) $\omega(x,y,z) := \exp(\Delta)$ relates the two parenthesizations via
		\begin{equation}
			\label{eq:assoziator_def}
			\text{Op}_{2; \, k_{xy}, k_z}^{(a)}\!\Big[ \text{Op}_{2; \, k_x, k_y}^{(a)}[x,y], \, z \Big] = \omega(x,y,z) \cdot \text{Op}_{2; \, k_x, k_{yz}}^{(a)}\!\Big[ x, \, \text{Op}_{2; \, k_y, k_z}^{(a)}[y,z] \Big],
		\end{equation}
		with exponent defect $\Delta = \frac{2\pi i}{a}\big(k_{xy}\Log_{k_z}(z) - k_{yz}\Log_{k_x}(x)\big)$, and enjoys the following properties:
		\begin{enumerate}
			\item \textbf{Centrality with respect to $y$ and modulus-phase factorization:} The associator $\omega(x,y,z) = \omega(x,z)$ is independent of the middle operand $y$. With total angles $\theta_x := \Arg(x) + 2\pi k_x$ and $\theta_z := \Arg(z) + 2\pi k_z$, it factors into a real amplitude scaling factor $|\omega|$ and a phase factor $e^{i\Phi}$:
			\begin{align}
				|\omega(x,z)| &= \exp\!\left( -\frac{2\pi}{a}\big( k_{xy}\theta_z - k_{yz}\theta_x \big) \right), \label{eq:omega_betrag} \\
				\Phi(x,z) &= \frac{2\pi}{a}\Big( k_{xy}\ln|z| - k_{yz}\ln|x| \Big) \pmod{2\pi}. \label{eq:omega_phase}
			\end{align}
			\item \textbf{Characterization of the associativity kernel:} Strict associativity ($\omega(x,z) = 1$) holds if and only if $\frac{1}{a}\big(k_{xy}\Log_{k_z}(z) - k_{yz}\Log_{k_x}(x)\big) \in \mathbb{Z}$. In particular, for symmetric intermediate indices $k_{xy} = k_{yz} = k \in \mathbb{Z}\setminus\{0\}$, the associativity kernel $\ker(\omega) := \{(x,z) \in S^2 \mid \omega(x,z) = 1\}$ is characterized by
			\begin{equation}
				\label{eq:invarianz_gitter}
				\Arg(z) = \Arg(x), \quad k_z = k_x \qquad \text{and} \qquad \frac{|z|}{|x|} = e^{\frac{a n}{k}} \quad (n \in \mathbb{Z}).
			\end{equation}
			Under the double-fixed-point calibration $a = \frac{\ln(2)}{2}$, this commensurability condition forms the dyadic lattice $|z|/|x| = 2^{\frac{n}{2k}}$.
		\end{enumerate}
	\end{theorem}
	
	\begin{proof}
		\textbf{Part 1:} Substituting $\Log_{k_{xy}}\!\big(\text{Op}_{2; \, k_x, k_y}^{(a)}[x,y]\big) = \frac{1}{a}\Log_{k_x}(x)\Log_{k_y}(y) + 2\pi i k_{xy}$ into the left-hand side of \eqref{eq:assoziator_def} and the analogous expression with $k_{yz}$ into the right-hand side gives the total exponents
		\begin{align*}
			E_1 &= \frac{1}{a^2}\Log_{k_x}(x)\Log_{k_y}(y)\Log_{k_z}(z) + \frac{2\pi i k_{xy}}{a}\Log_{k_z}(z), \\
			E_2 &= \frac{1}{a^2}\Log_{k_x}(x)\Log_{k_y}(y)\Log_{k_z}(z) + \frac{2\pi i k_{yz}}{a}\Log_{k_x}(x).
		\end{align*}
		Taking the difference $\Delta := E_1 - E_2$ cancels the trilinear leading term together with the middle operand $y$. Substituting $\Log_{k_z}(z) = \ln|z| + i\theta_z$ and $\Log_{k_x}(x) = \ln|x| + i\theta_x$ and separating real and imaginary parts gives
		\begin{equation*}
			\Delta = -\frac{2\pi}{a}\big(k_{xy}\theta_z - k_{yz}\theta_x\big) + i\,\frac{2\pi}{a}\big(k_{xy}\ln|z| - k_{yz}\ln|x|\big),
		\end{equation*}
		from which \eqref{eq:omega_betrag} and \eqref{eq:omega_phase} follow immediately via $\omega = \exp(\operatorname{Re}\Delta)\exp(i\operatorname{Im}\Delta)$. (If every intermediate result on $S$ is evaluated via the principal value $\Log_0$, then $k_{xy} = -\lfloor(\operatorname{Im}(\frac{1}{a}\Log_{k_x}(x)\Log_{k_y}(y)) + \pi)/(2\pi)\rfloor$ is intrinsically determined by exiting the fundamental strip $(-\pi, \pi]$.)
		
		\textbf{Part 2:} The condition $\omega(x,z) = 1 \iff \Delta \in 2\pi i \mathbb{Z}$ simultaneously requires $|\omega| = 1$ and $\Phi \in 2\pi\mathbb{Z}$. For $k_{xy} = k_{yz} = k \neq 0$, $|\omega| = 1$ in \eqref{eq:omega_betrag} yields $\theta_z = \theta_x$, which, because $\Arg(x), \Arg(z) \in (-\pi, \pi]$, is equivalent to $\Arg(z) = \Arg(x)$ and $k_z = k_x$. Finally, $\Phi = \frac{2\pi k}{a}\ln(|z|/|x|) = 2\pi n$ ($n \in \mathbb{Z}$) gives $\ln(|z|/|x|) = \frac{an}{k}$, establishing \eqref{eq:invarianz_gitter}.
	\end{proof}
	
	\section{Cumulant Expansion and the General Covariance Theorem for Commutator Defects}
	
	\subsection{Cumulant Expansion and Sign Reversal of the Quadratic Form}
	
	By virtue of its structure as a cumulant-generating function \cite{chapman1970, grad1949, jarzynski1997}, the Taylor expansion of the mean $S_{\text{Op}_{-1}}(x_1, \dots, x_n) = \frac{1}{a}\ln\!\big(\frac{1}{n}\sum_{i=1}^n e^{a x_i}\big)$ around the arithmetic mean $\mu = S_{\text{Op}_0}(\mathbf{x}) = \frac{1}{n}\sum_{i=1}^n x_i$ in powers of $a$ yields the empirical cumulants of the sample:
	\begin{equation}
		\label{eq:kumulanten_ne}
		S_{\text{Op}_{-1}}(x_1, \dots, x_n) = \mu + \frac{a}{2}\!\left( \frac{1}{n}\sum_{i=1}^n (x_i - \mu)^2 \right) + \frac{a^2}{6}\!\left( \frac{1}{n}\sum_{i=1}^n (x_i - \mu)^3 \right) + \mathcal{O}(a^3).
	\end{equation}
	For $n = 2$ arguments $(x_1, x_2)$ with midpoint $m = \frac{x_1+x_2}{2}$ and half-range $s = \frac{x_1-x_2}{2}$, we have the exact identity $S_{\text{Op}_{-1}}(x_1, x_2) = m + \frac{1}{a}\ln\cosh(as) = m + \frac{a}{2}s^2 + \mathcal{O}(s^4) > S_{\text{Op}_0}(x_1, x_2)$ ($s \neq 0$), while for independent, identically distributed random variables $X_i$ with $\mathbb{E}[e^{a X_i}] = \mu_Y < \infty$ and $\Var(e^{a X_i}) = \sigma_Y^2 < \infty$, the Central Limit Theorem \cite{feller1968} and the delta method \cite{cramer1946} yield asymptotic normality $\sqrt{n}\big(S_{\text{Op}_{-1}}(X_1, \dots, X_n) - \frac{1}{a}\ln\mu_Y\big) \xrightarrow{d} \mathcal{N}\big(0, \frac{\sigma_Y^2}{a^2\mu_Y^2}\big)$. By contrast, in logarithmic coordinates $\eta_i = \ln x_i$, $\bar{\eta} = \frac{\eta_1+\eta_2}{2}$, $\delta = \frac{\eta_1-\eta_2}{2}$, the operation $\text{Op}_2^{(a)}[x_1, x_2]$ leads to the indefinite bilinear form $\ln\!\big(\text{Op}_2^{(a)}[x_1, x_2]\big) = \frac{1}{a}\eta_1\eta_2 = \frac{1}{a}(\bar{\eta}^2 - \delta^2)$, where any dispersion $\delta \neq 0$ reduces the output value.
	
	\subsection{The General Covariance Pushforward Theorem for the Commutator Defect}
	
	Because the mean $S_{\text{Op}_m}$ is linearized with respect to $\text{Op}_m^{(a)}$, it does not in general commute with the superordinate operation $\text{Op}_{m+1}^{(a)}$. We define this commutator defect for general sample sizes $n \in \mathbb{N}$ and then prove as our main theorem that the defects across all levels $m \in \mathbb{Z}$ are mapped onto the classical empirical covariance by the isomorphism $f_a^{(m)}$.
	
	\begin{definition}[Level-specific commutator defect]
		\label{def:spread}
		Let $\mathbf{x} = (x_1, \dots, x_n) \in D_m^n$ and $\mathbf{y} = (y_1, \dots, y_n) \in D_m^n$ be two sample vectors at level $m \in \mathbb{Z}$. The \textbf{commutator defect} (hierarchical spread) $\text{Sp}_{\text{Op}_{m+1}}(\mathbf{x}, \mathbf{y})$ between the mean $S_{\text{Op}_m}$ and the superordinate operation $\text{Op}_{m+1}^{(a)}$ is defined implicitly by
		\begin{equation}
			\label{eq:spread_allgemein}
			\begin{aligned}
				\text{Op}_m^{(a)}\!\Big[ \text{Sp}_{\text{Op}_{m+1}}(\mathbf{x}, \mathbf{y}), \; &\text{Op}_{m+1}^{(a)}\big[ S_{\text{Op}_m}(\mathbf{x}), \, S_{\text{Op}_m}(\mathbf{y}) \big] \Big] \\
				&= S_{\text{Op}_m}\!\Big( \text{Op}_{m+1}^{(a)}[x_1, y_1], \dots, \text{Op}_{m+1}^{(a)}[x_n, y_n] \Big).
			\end{aligned}
		\end{equation}
	\end{definition}
	
	\begin{theorem}[General covariance pushforward theorem for $\text{Sp}_{\text{Op}_{m+1}}$]
		\label{thm:kovarianz_pushforward}
		Let $a > 0$, $m \in \mathbb{Z}$, and $n \in \mathbb{N}$. For arbitrary sample vectors $\mathbf{x}, \mathbf{y} \in D_m^n$, let $u_i := f_a^{(m)}(x_i) \in \mathbb{R}$ and $v_i := f_a^{(m)}(y_i) \in \mathbb{R}$ denote the coordinates transformed into the linear reference space, and let
		\begin{equation*}
			\Cov(\mathbf{u}, \mathbf{v}) := \frac{1}{n}\sum_{i=1}^n u_i v_i - \left(\frac{1}{n}\sum_{i=1}^n u_i\right)\!\left(\frac{1}{n}\sum_{i=1}^n v_i\right) = \frac{1}{n}\sum_{i=1}^n (u_i - \bar{u})(v_i - \bar{v})
		\end{equation*}
		be their empirical sample covariance. Then at every level $m \in \mathbb{Z}$ the exact identity
		\begin{equation}
			\label{eq:master_kovarianz_id}
			f_a^{(m)}\!\Big( \text{Sp}_{\text{Op}_{m+1}}(\mathbf{x}, \mathbf{y}) \Big) = \Cov\!\Big( f_a^{(m)}(\mathbf{x}), \, f_a^{(m)}(\mathbf{y}) \Big)
		\end{equation}
		holds. Whenever $\Cov(f_a^{(m)}(\mathbf{x}), f_a^{(m)}(\mathbf{y}))$ lies in the range of $f_a^{(m)}$ (for $m \ge 0$ for arbitrary correlations in $\mathbb{R}$, and for $m = -1$ under positive correlation $\Cov(e^{a\mathbf{x}}, e^{a\mathbf{y}}) > 0$), the commutator defect admits the closed explicit representation
		\begin{equation}
			\label{eq:master_spread_explizit}
			\text{Sp}_{\text{Op}_{m+1}}(\mathbf{x}, \mathbf{y}) = f_a^{(-m)}\!\Big( \Cov\!\big( f_a^{(m)}(\mathbf{x}), \, f_a^{(m)}(\mathbf{y}) \big) \Big).
		\end{equation}
		In particular, $\text{Sp}_{\text{Op}_{m+1}}(\mathbf{x}, \mathbf{y})$ equals the neutral element $n_{\text{Op}_m} = f_a^{(-m)}(0)$ of the subordinate operation if and only if the transformed samples $f_a^{(m)}(\mathbf{x})$ and $f_a^{(m)}(\mathbf{y})$ are uncorrelated ($\Cov(f_a^{(m)}(\mathbf{x}), f_a^{(m)}(\mathbf{y})) = 0$).
	\end{theorem}
	
	\begin{proof}
		By Definition~\ref{def:isomorphismus_allg} and Theorem~\ref{thm:master_normierung}, for each $m \in \mathbb{Z}$ the bijection $f_a^{(m)}$ simultaneously maps the three operations involved onto the standard operations on $\mathbb{R}$:
		\begin{align*}
			f_a^{(m)}\!\big(\text{Op}_m^{(a)}[x, y]\big) &= f_a^{(m)}(x) + f_a^{(m)}(y), \\
			f_a^{(m)}\!\big(\text{Op}_{m+1}^{(a)}[x, y]\big) &= f_a^{-1}\!\Big( f_a^{(m+1)}(x) + f_a^{(m+1)}(y) \Big) = f_a^{(m)}(x) \cdot f_a^{(m)}(y), \\
			f_a^{(m)}\!\big(S_{\text{Op}_m}(z_1, \dots, z_n)\big) &= \frac{1}{n}\sum_{i=1}^n f_a^{(m)}(z_i).
		\end{align*}
		Applying $f_a^{(m)}$ to both sides of the defining equation \eqref{eq:spread_allgemein} and substituting $u_i = f_a^{(m)}(x_i)$ and $v_i = f_a^{(m)}(y_i)$ transforms the left-hand side into
		\begin{equation*}
			f_a^{(m)}\!\Big(\text{Sp}_{\text{Op}_{m+1}}(\mathbf{x}, \mathbf{y})\Big) + \left(\frac{1}{n}\sum_{i=1}^n u_i\right)\!\left(\frac{1}{n}\sum_{i=1}^n v_i\right),
		\end{equation*}
		while the right-hand side of \eqref{eq:spread_allgemein} yields the mean of the products:
		\begin{align*}
			f_a^{(m)}\!\left( S_{\text{Op}_m}\!\Big(\text{Op}_{m+1}^{(a)}[x_1, y_1], \dots, \text{Op}_{m+1}^{(a)}[x_n, y_n]\Big) \right) &= \frac{1}{n}\sum_{i=1}^n f_a^{(m)}\!\big(\text{Op}_{m+1}^{(a)}[x_i, y_i]\big) \\
			&= \frac{1}{n}\sum_{i=1}^n u_i v_i.
		\end{align*}
		Solving for $f_a^{(m)}\big(\text{Sp}_{\text{Op}_{m+1}}(\mathbf{x}, \mathbf{y})\big)$ immediately proves \eqref{eq:master_kovarianz_id}, and applying $f_a^{(-m)}$ yields \eqref{eq:master_spread_explizit}.
	\end{proof}
	
	Theorem~\ref{thm:kovarianz_pushforward} acts as a structural pushforward (or translation) theorem under the level isomorphism $f_a^{(m)}$: it shows that the discrepancy measures appearing across different scales for $m \in \{-1, 0, 1\}$ are the transformed images of one and the same bilinear covariance form (for $n = 2$, $\Cov(\mathbf{u}, \mathbf{v}) = \frac{1}{4}(u_1 - u_2)(v_1 - v_2)$ in each case).
	
	\paragraph{Illustrative Applications:}
	\begin{enumerate}
		\item \textbf{Level $m = -1$ (Logarithmic level bias, $\text{Op}_0 = Ad, \, \text{Op}_{-1} = Ne$):} \\
		With $f_a^{(-1)}(x) = e^{ax}$ and $f_a^{(1)}(u) = \frac{1}{a}\ln(u)$, equation \eqref{eq:master_spread_explizit} directly gives $\text{Sp}_{\text{Op}_0}(\mathbf{x}, \mathbf{y}) = \frac{1}{a}\ln\Cov(e^{a\mathbf{x}}, e^{a\mathbf{y}})$. In acoustic decibel scaling ($a_{\text{acoust}} = \frac{\ln 10}{10}$), superimposing a constant and an intermittent level series yields the energetic level bias $\Delta L_{\text{bias}}$ of technical acoustics \cite{moeser2012, din9612}.
		\item \textbf{Level $m = 0$ (Classical sample covariance, $\text{Op}_1 = Mu, \, \text{Op}_0 = Ad$):} \\
		Since $f_a^{(0)}(x) = x$ is the identity map, equation \eqref{eq:master_spread_explizit} reduces to $\text{Sp}_{\text{Op}_1}(\mathbf{x}, \mathbf{y}) = \Cov(\mathbf{x}, \mathbf{y})$ (specifically $\frac{1}{4}(x_1-x_2)(y_1-y_2)$ for $n=2$) \cite{feller1968}.
		\item \textbf{Level $m = 1$ (Exponential logarithmic covariance factor, $\text{Op}_2 = At, \, \text{Op}_1 = Mu$):} \\
		With $f_a^{(1)}(x) = \frac{1}{a}\ln(x)$ and $f_a^{(-1)}(u) = e^{au}$, equation \eqref{eq:master_spread_explizit} gives for arbitrary $n \in \mathbb{N}$:
		\begin{equation}
			\label{eq:spread_at_ito}
			\text{Sp}_{\text{Op}_2}(\mathbf{x}, \mathbf{y}) = \exp\!\Big( a \cdot \Cov\!\big(\tfrac{1}{a}\ln\mathbf{x}, \, \tfrac{1}{a}\ln\mathbf{y}\big) \Big) = \exp\!\left( \frac{1}{a}\,\Cov(\ln\mathbf{x}, \, \ln\mathbf{y}) \right).
		\end{equation}
		For identical time series $\mathbf{x} = \mathbf{y}$ with logarithmic sample variance $\sigma_{\ln}^2 := \Var(\ln\mathbf{x})$, we obtain $\text{Sp}_{\text{Op}_2}(\mathbf{x}, \mathbf{x}) = \exp(\frac{1}{a}\sigma_{\ln}^2)$, capturing the exponential variance structure of log-normal growth processes in portfolio theory \cite{luenberger1998, markowitz1952}.
	\end{enumerate}
	
	\begin{table}[htbp]
		\centering
		\small
		\renewcommand{\arraystretch}{1.28}
		\caption{Systematic overview of the levels $U = \text{Op}_m^{(a)}$ (with adjacent level $O = \text{Op}_{m+1}^{(a)}$ and calibration $a = \frac{\ln(2)}{2}$)}
		\label{tab:hierarchie_uebersicht}
		\vspace{1ex}
		\begin{tabular*}{\textwidth}{@{\extracolsep{\fill}}lll@{}}
			\toprule
			\textbf{Criterion (Part A: Lower levels)} & \textbf{Level $m = -1$} & \textbf{Level $m = 0$} \\
			\midrule
			Operation $U = \text{Op}_m^{(a)}$ & $\frac{1}{a}\ln(e^{ax}+e^{ay})$ & $x + y$ \\
			Transformed coordinate $f_a^{(m)}(x)$ & $e^{ax}$ & $x$ \\
			Real group/monoid space & $\mathbb{R} \cup \{-\infty\}$ & $\mathbb{R}$ \\
			Neutral element $n_{\text{Op}_m}$ & $-\infty$ & $0$ \\
			Absorbing element & $+\infty$ \textit{(in $\overline{\mathbb{R}}$)} & $-\infty$ \\
			Mean $S_{\text{Op}_m}(\mathbf{x})$ & $\frac{1}{a}\ln(\frac{1}{n}\sum e^{ax_i})$ & $\frac{1}{n}\sum x_i$ \\
			Normalizer $N_{\text{Op}_{m+1}}(n)$ & $-\frac{1}{a}\ln(n)$ & $\frac{1}{n}$ \\
			Commutator defect $\text{Sp}_{\text{Op}_{m+1}}$ & $\frac{1}{a}\ln\Cov(e^{a\mathbf{x}}, e^{a\mathbf{y}})$ & $\Cov(\mathbf{x}, \mathbf{y})$ \\
			\bottomrule
		\end{tabular*}
		
		\vspace{2ex}
		\begin{tabular*}{\textwidth}{@{\extracolsep{\fill}}lll@{}}
			\toprule
			\textbf{Criterion (Part B: Upper levels)} & \textbf{Level $m = +1$} & \textbf{Level $m = +2$} \\
			\midrule
			Operation $U = \text{Op}_m^{(a)}$ & $x \cdot y$ & $\exp(\frac{1}{a}\ln x \ln y)$ \\
			Transformed coordinate $f_a^{(m)}(x)$ & $\frac{1}{a}\ln(x)$ & $\frac{1}{a}\ln(\frac{1}{a}\ln x)$ \\
			Real group/monoid space & $\mathbb{R}_{>0}$ & $\mathbb{R}_{>0}\setminus\{1\}$ \\
			Neutral element $n_{\text{Op}_m}$ & $1$ & $e^a = \sqrt{2}$ \\
			Absorbing element & $0$ & $1$ \\
			Mean $S_{\text{Op}_m}(\mathbf{x})$ & $(\prod x_i)^{1/n}$ & $\exp\!\big(a\prod (\frac{\ln x_i}{a})^{1/n}\big)$ \\
			Normalizer $N_{\text{Op}_{m+1}}(n)$ & $e^{a/n} = 2^{\frac{1}{2n}}$ & $e^{a e^{a/n}} = 2^{\frac{1}{2}\cdot 2^{\frac{1}{2n}}}$ \\
			Commutator defect $\text{Sp}_{\text{Op}_{m+1}}$ & $\exp\!\big(\frac{1}{a}\Cov(\ln\mathbf{x}, \ln\mathbf{y})\big)$ & \begin{tabular}[t]{@{}l@{}}$f_a^{(-2)}\!\big(\Cov(f_a^{(2)}(\mathbf{x}),$ \\ $\qquad f_a^{(2)}(\mathbf{y}))\big)$\end{tabular} \\
			\bottomrule
		\end{tabular*}
	\end{table}
	
	\section{Summary and Outlook}
	
	The present investigation has shown that the hierarchy of operations $\text{Op}_m^{(a)}$ generated by iterated logarithmic conjugation exhibits three specific structural properties beyond the mere transport of abelian group structures:
	\begin{enumerate}
		\item Within continuously conjugated group operations, the iterated logarithms $\varphi_m = f_a^{(m)}$ are, by Theorem~\ref{thm:charakterisierung_mittel}, precisely those generators whose quasi-arithmetic means $S_{\text{Op}_m}$ are internally normalized by the superordinate adjacent operation $\text{Op}_{m+1}^{(a)}$ via the universal factor $N_{\text{Op}_{m+1}}(n) = f_a^{(-m)}(1/n)$ and are simultaneously bi-invariant with respect to $\text{Op}_m^{(a)}$ and $\text{Op}_{m+1}^{(a)}$.
		\item The real existence of these normalization factors terminates for $n \ge 2$ below $m = -1$ upon leaving the positive real domain $\mathbb{R}_{>0}$ (Theorem~\ref{thm:master_normierung}, Part~2), motivating the analytic continuation onto the Riemann surface of the logarithm, where $\text{Op}_{-1}^{(a)}$ remains inner-branch invariant under nested evaluation (Theorem~\ref{thm:komplex_ne_eigenschaften}), whereas $\text{Op}_2^{(a)}$ generates an associator $\omega(x,z)$ factoring into amplitude scaling and phase rotation (Theorem~\ref{thm:komplex_assoziativ}).
		\item By the covariance pushforward theorem (Theorem~\ref{thm:kovarianz_pushforward}), the level-specific commutator defects $\text{Sp}_{\text{Op}_{m+1}}(\mathbf{x}, \mathbf{y})$ admit the closed representation $f_a^{(-m)}\!\big(\Cov(f_a^{(m)}(\mathbf{x}), f_a^{(m)}(\mathbf{y}))\big)$ for all $m \in \mathbb{Z}$ and $n \in \mathbb{N}$.
	\end{enumerate}
	In future work, this formalism will be applied to stochastic bifurcation models, in particular to the transition from unimodal binomial distributions to bimodal arcsine distributions at constant arithmetic mean $S_{\text{Op}_0}$ and monotonically increasing defect $S_{\text{Op}_{-1}} - S_{\text{Op}_0}$.
	
	\section*{Acknowledgments}
	I warmly thank P.~Collignon, M.~Girbert, S.~Kraußhar, as well as A.~and D.~Legatiuk for valuable and inspiring discussions, and I dedicate this article to my sons Tobias and Haibo.

\end{document}